\input amstex\documentstyle{amsppt}  
\pagewidth{12.5cm}\pageheight{19cm}\magnification\magstep1
\topmatter
\title Action of longest element on a Hecke algebra cell module\endtitle
\author G. Lusztig\endauthor
\address{Department of Mathematics, M.I.T., Cambridge, MA 02139}\endaddress
\dedicatory{Dedicated to the memory of Robert Steinberg}\enddedicatory
\thanks{Supported in part by National Science Foundation grant DMS-1303060.}\endthanks
\endtopmatter   
\document

\define\Irr{\text{\rm Irr}}

\define\mpb{\medpagebreak}

\define\da{\dagger}

\define\si{\sim}

\define\lb{\linebreak}

\define\op{\oplus}
   
\define\part{\partial}
\define\emp{\emptyset}

\define\n{\notin}

\define\m{\mapsto}

\define\lra{\leftrightarrow}

\define\sub{\subset}

\define\nl{\newline}
\redefine\i{^{-1}}

\define\ov{\overline}
\define\ot{\otimes}

\define\tr{\text{\rm tr}}

\redefine\spa{\spadesuit}

\define\g{\gamma}
\redefine\d{\delta}
\define\e{\epsilon}

\define\io{\iota}

\define\ph{\phi}

\define\s{\sigma}

\define\th{\theta}

\redefine\l{\lambda}

\define\x{\xi}

\define\vt{\vartheta}

\redefine\G{\Gamma}

\redefine\aa{\bold a}

\define\boc{\bold c}

\define\CC{\bold C}

\define\NN{\bold N}

\define\ZZ{\bold Z}

\define\ca{\Cal A}

\define\cd{\Cal D}

\define\ch{\Cal H}

\define\cl{\Cal L}

\define\car{\Cal R}

\define\fT{\frak T}

\define\sha{\sharp}

\define\BFO{BFO}
\define\EW{EW}
\define\KL{KL}
\define\ORA{L1}
\define\CBII{L2}
\define\HEC{L3}
\define\CDGVII{L4}
\define\MAT{MA}

\head Introduction\endhead
\subhead 0.1\endsubhead
The Hecke algebra $\ch$ (over $\ca=\ZZ[v,v\i]$, $v$ an indeterminate) of a finite Coxeter group $W$ has two 
bases as an $\ca$-module: the standard basis $\{T_x;x\in W\}$ and the basis $\{C_x;x\in W\}$ introduced in 
\cite{\KL}. The second basis determines a decomposition of $W$ into two-sided cells and a partial order for 
the set of two-sided cells, see \cite{\KL}. Let $l:W@>>>\NN$ be the length function, let $w_0$ be the 
longest element of $W$ and let $\boc$ be a two-sided cell. Let $a$ (resp. $a'$) be the value of the 
$\aa$-function \cite{\HEC, 13.4} on $\boc$ (resp. on $w_0\boc$). The following result was proved by Mathas in \cite{\MAT}.

(a) {\it There exists a unique permutation $u\m u^*$ of $\boc$ such that for any $u\in\boc$ we have 
$T_{w_0}(-1)^{l(u)}C_u=(-1)^{l(w_0)+a'}v^{-a+a'}(-1)^{l(u^*)}C_{u^*}$ plus an $\ca$-linear combination of 
elements $C_{u'}$ with $u'$ in a two-sided cell strictly smaller than $\boc$. Moreover, for any $u\in\boc$ 
we have $(u^*)^*=u$.}
\nl
A related (but weaker) result appears in \cite{\ORA, (5.12.2)}.

A result similar to (a) which concerns canonical bases in representations of quantum groups appears in 
\cite{\CBII, Cor. 5.9}; now, in the case where $W$ is of type $A$, (a) can be deduced from {\it loc.cit.} 
using the fact that irreducible representations of the Hecke algebra of type $A$ (with their canonical 
bases) can be realized as $0$-weight spaces of certain irreducible representations of a quantum group with
their canonical bases.

As R. Bezrukavnikov pointed out to the author, (a) specialized for $v=1$  (in the group algebra of $W$ 
instead of $\ch$) and assuming that $W$ is crystallographic can be deduced from \cite{\BFO, Prop. 4.1} (a 
statement about Harish-Chandra modules), although it is not explicitly stated there.

In this paper we shall prove a generalization of (a) which applies to the Hecke algebra associated to $W$ and
any weight function assumed to satisfy the properties P1-P15 in \cite{\HEC,\S14}, see Theorem 2.3; (a)
corresponds to the special case where the weight function is equal to the length function. As an 
application we show that the image of $T_{w_0}$ in the asymptotic Hecke algebra is given by a simple formula
(see 2.8).

I thank Matthew Douglass for bringing the paper \cite{\MAT} to my attention. I thank the referee for
helpful comments.

\subhead 0.2\endsubhead
{\it Notation. } $W$ is a finite Coxeter group; the set of simple reflections is denoted by $S$. We shall 
adopt many notations of \cite{\HEC}. Let $\le$ be the standard partial order on $W$. Let $l:W@>>>\NN$ be the
 length function of $W$ and let $L:W@>>>\NN$ be a weight function (see \cite{\HEC, 3.1}) 
that is, a function such that $L(ww')=L(w)+L(w')$ for any $w,w'$ in $W$ such that $l(ww')=l(w)+l(w')$; we 
assume that $L(s)>0$ for any $s\in S$. Let $w_0,\ca$ be as in 0.1 and let $\ch$ be the 
Hecke algebra over $\ca$ associated to $W,L$ as in \cite{\HEC, 3.2}; we shall assume that properties
P1-P15 in \cite{\HEC, \S14} are satisfied. (This holds automatically if $L=l$ by \cite{\HEC,\S15} using the 
results of \cite{\EW}. This also holds in the quasisplit case, see 
\cite{\HEC,\S16}.)
We have $\ca\sub\ca'\sub K$ where $\ca'=\CC[v,v\i],K=\CC(v)$. Let $\ch_K=K\ot_\ca\ch$ (a $K$-algebra).
Recall that $\ch$ has an $\ca$-basis $\{T_x;x\in W\}$, see \cite{\HEC, 3.2} and an $\ca$-basis 
$\{c_x;x\in W\}$, see \cite{\HEC, 5.2}. 
For $x\in W$ we have $c_x=\sum_{y\in W}p_{y,x}T_y$ and $T_x=\sum_{y\in W}(-1)^{l(xy)}p_{w_0x,w_0y}c_y$
(see \cite{\HEC, 11.4}) where $p_{x,x}=1$ and $p_{y,x}\in v\i\ZZ[v\i]$ for $y\ne x$.
We define preorders $\le_\cl,\le_\car,\le_{\cl\car}$ on $W$ in terms
of $\{c_x;x\in W\}$ as in \cite{\HEC, 8.1}. Let $\si_\cl,\si_\car,\si_{\cl\car}$ be the corresponding 
equivalence relations on $W$, see \cite{\HEC, 8.1} (the equivalence classes are called left cells, right
cells, two-sided cells). Let $\bar{}:\ca@>>>\ca$ be the ring involution such that $\ov{v^n}=v^{-n}$ for 
$n\in\ZZ$.
Let $\bar{}:\ch@>>>\ch$ be the ring involution such that $\ov{fT_x}=\bar fT_{x\i}\i$ for $x\in W,f\in\ca$.
For $x\in W$ we have $\ov{c_x}=c_x$.
Let $h\m h^\da$ be the algebra automorphism of $\ch$ or of $\ch_K$ given by $T_x\m(-1)^{l(x)}T_{x\i}\i$ for 
all $x\in W$, see \cite{\HEC, 3.5}.
Then the basis $\{c_x^\da;x\in W\}$ of $\ch$ is defined. (In the case where $L=l$, for any $x$ we have 
$c_x^\da=(-1)^{l(x)}C_x$ where $C_x$ is as in 0.1.)
Let $h\m h^\flat$ be the algebra antiautomorphism of $\ch$ given by $T_x\m T_{x\i}$ for all $x\in W$, see 
\cite{\HEC, 3.5}; for $x\in W$ we have $c_x^\flat=c_{x\i}$, see \cite{\HEC, 5.8}.
For $x,y\in W$ we have $c_xc_y=\sum_{z\in W}h_{x,y,z}c_z$, $c_x^\da c_y^\da=\sum_{z\in W}h_{x,y,z}c_z^\da$,
where $h_{x,y,z}\in\ca$.  
For any $z\in W$ there is a unique number $\aa(z)\in\NN$ such that for any $x,y$ in $W$ we have
$$h_{x,y,z}=\g_{x,y,z\i}v^{\aa(z)}+\text{strictly smaller powers of } v$$
where $g_{x,y,z\i}\in\ZZ$ and $g_{x,y,z\i}\ne0$ for some $x,y$ in $W$. We have also
$$h_{x,y,z}=\g_{x,y,z\i}v^{-\aa(z)}+\text{strictly larger powers of } v.$$
Moreover $z\m\aa(z)$ is constant on
any two-sided cell. The free abelian group $J$ with basis $\{t_w;w\in W\}$ has an associative ring structure
given by $t_xt_y=\sum_{z\in W}\g_{x,y,z\i}t_z$; it has a unit element of the form
$\sum_{d\in\cd}n_dt_d$ where $\cd$ is a subset of $W$ consisting of certain elements with square $1$
and $n_d=\pm1$. Moreover for $d\in\cd$ we have $n_d=\g_{d,d,d}$.

For any $x\in W$ there is a unique element $d_x\in\cd$ such that $x\si_\cl d_x$.
For a commutative ring $R$ with $1$ we set $J_R=R\ot J$ (an $R$-algebra).

There is a unique $\ca$-algebra homomorphism $\ph:\ch@>>>J_\ca$ such that
$\ph(c_x^\da)=\sum_{d\in\cd,z\in W;d_z=d}h_{x,d,z}n_dt_z$ for any $x\in W$.
After applying $\CC\ot_\ca$ to $\ph$ (we regard $\CC$ as an $\ca$-algebra via $v\m1$), $\ph$ becomes a
$\CC$-algebra isomorphism $\ph_\CC:\CC[W]@>\si>>J_\CC$ (see \cite{\HEC, 20.1(e)}). After applying $K\ot_\ca$
to $\ph$, $\ph$ becomes a $K$-algebra isomorphism $\ph_K:\ch_K@>\si>>J_K$ (see \cite{\HEC, 20.1(d)}).

For any two-sided cell $\boc$ let $\ch^{\le\boc}$ (resp. $\ch^{<\boc}$) be the $\ca$-submodule of $\ch$ 
spanned by $\{c_x^\da,x\in W,x\le_{\cl\car}x'\text{ for some }x'\in\boc\}$
(resp. $\{c_x^\da,x\in W,x<_{\cl\car}x'\text{ for some }x'\in\boc\}$). Note that 
$\ch^{\le\boc},\ch^{<\boc}$ are two-sided ideals in $\ch$. Hence $\ch^\boc:=\ch^{\le\boc}/\ch^{<\boc}$ is an
$\ch,\ch$ bimodule. It has an $\ca$-basis $\{c_x^\da,x\in\boc\}$. Let $J^\boc$ be the subgroup of $J$ 
spanned by $\{t_x;x\in\boc\}$. This is a two-sided ideal of $J$.
Similarly, $J^\boc_\CC:=\CC\ot J^\boc$ is a two-sided ideal of $J_\CC$ and
$J^\boc_K:=K\ot J^\boc$ is a two-sided ideal of $J_K$.

We write $E\in\Irr W$ whenever $E$ is a simple $\CC[W]$-module. We can view $E$ as a (simple) $J_\CC$-module
$E_\spa$ via the isomorphism $\ph_\CC\i$. Then the (simple) $J_K$-module $K\ot_\CC E_\spa$ can be viewed 
as a (simple) $\ch_K$-module $E_v$ via the isomorphism $\ph_K$.
Let $E^\da$ be the simple $\CC[W]$-module which coincides with $E$ as a $\CC$-vector space
but with the $w$ action on $E^\da$ (for $w\in W$) being $(-1)^{l(w)}$ times the $w$-action on $E$.
Let $\aa_E\in\NN$ be as in \cite{\HEC, 20.6(a)}.

\head 1. Preliminaries\endhead
\subhead 1.1\endsubhead
Let $\s:W@>>>W$ be the automorphism given by $w\m w_0ww_0$; it satisfies $\s(S)=S$ and it extends to a 
$\CC$-algebra isomorphism $\s:\CC[W]@>>>\CC[W]$. For $s\in S$ we have 
$l(w_0)=l(w_0s)+l(s)=l(\s(s))+l(\s(s)w_0)$ hence $L(w_0)=L(w_0s)+L(s)=L(\s(s))+L(\s(s)w_0)=L(\s(s))+L(w_0s)$
so that $L(\s(s))=L(s)$.
It follows that $L(\s(w))=L(w)$ for all $w\in W$ and that we have an $\ca$-algebra automorphism
$\s:\ch@>>>\ch$ where $\s(T_w)=T_{\s(w)}$ for any $w\in W$. This extends to a
$K$-algebra isomorphism $\s:\ch_K@>>>\ch_K$. We have $\s(c_w)=c_{\s(w)}$ for any $w\in W$. For any 
$h\in\ch$ we have $\s(h^\da)=(\s(h))^\da$. Hence we have $\s(c_w^\da)=c_{\s(w)}^\da$ for any $w\in W$.
We have $h_{\s(x),\s(y),\s(z)}=h_{x,y,z}$ for all
$x,y,z\in W$. It follows that $\aa(\s(w))=\aa(w)$ for all $w\in W$ and 
$\g_{\s(x),\s(y),\s(z)}=\g_{x,y,z}$ for all $x,y,z\in W$ so that we have a ring isomorphism
$\s:J@>>>J$ where $\s(t_w)=t_{\s(w)}$ for any $w\in W$. This extends to
an $\ca$-algebra isomorphism $\s:J_\ca@>>>J_\ca$, to a $\CC$-algebra isomorphism $\s:J_\CC@>>>J_\CC$ and to 
a $K$-algebra isomorphism $\s:J_K@>>>J_K$.
From the definitions we see that $\ph:\ch@>>>J_\ca$ (see 0.2) satisfies $\ph\s=\s\ph$.
Hence $\ph_\CC$ satisfies $\ph_\CC\s=\s\ph_\CC$ and $\ph_K$ satisfies $\ph_K\s=\s\ph_K$.

We show:

(a) For $h\in\ch$ we have $\s(h)=T_{w_0}hT_{w_0}\i$. 
\nl
It is enough to show this for $h$ running through a set of algebra generators of $\ch$. Thus we can assume
that $h=T_s\i$ with $s\in S$. We must show that $T_{\s(s)}\i T_{w_0}=T_{w_0}T_s\i$: both sides are equal
to $T_{\s(s)w_0}=T_{w_0s}$.

\proclaim{Lemma 1.2} For any $x\in W$ we have $\s(x)\si_{\cl\car}x$.
\endproclaim
From 1.1(a) we deduce that $T_{w_0}c_xT_{w_0}\i=c_{\s(x)}$. In particular, $\s(x)\le_{\cl\car}x$. Replacing 
$x$ by $\s(x)$ we obtain $x\le_{\cl\car}\s(x)$. The lemma follows.

\subhead 1.3\endsubhead
Let $E\in\Irr W$. We define $\s_E:E@>>>E$ by $\s_E(e)=w_0e$ for $e\in E$. We have $\s_E^2=1$. For 
$e\in E,w\in W$, we have $\s_E(we)=\s(w)\s_E(e)$. We can view $\s_E$ as a vector space isomorphism 
$E_\spa@>\si>>E_\spa$. For $e\in E_\spa,w\in W$ we have $\s_E(t_we)=t_{\s(w)}\s_E(e)$. Now 
$\s_E:E_\spa@>>>E_\spa$ defines by extension of scalars a vector space isomorphism $E_v@>>>E_v$ denoted 
again by $\s_E$. It satisfies $\s_E^2=1$. For $e\in E_v,w\in W$ we have $\s_E(T_we)=T_{\s(w)}\s_E(e)$.

\proclaim{Lemma 1.4} Let $E\in\Irr W$. There is a unique (up to multiplication by a scalar in $K-\{0\}$) 
vector space isomorphism $g:E_v@>>>E_v$ such that $g(T_we)=T_{\s(w)}g(e)$ for all $w\in W,e\in E_v$. We can 
take for example $g=T_{w_0}:E_v@>>>E_v$ or $g=\s_E:E_v@>>>E_v$. Hence $T_{w_0}=\l_E\s_E:E_v@>>>E_v$ where 
$\l_E\in K-\{0\}$.
\endproclaim
The existence of $g$ is clear from the second sentence of the lemma. If $g'$ is another isomorphism
$g':E_v@>>>E_v$ such that $g'(T_we)=T_{\s(w)}g'(e)$ for all $w\in W,e\in E_v$ then for any $e\in E_v$ we 
have $g\i g'(T_we)=g\i T_{\s(w)}g'(e)=T_wg\i g'(e)$ and using Schur's lemma we see that $g\i g'$ is a 
scalar. This proves the first sentence of the lemma hence the third sentence of the lemma.

\subhead 1.5\endsubhead
Let $E\in\Irr W$. We have 
$$\sum_{x\in W}\tr(T_x,E_v)\tr(T_{x\i},E_v)=f_{E_v}\dim(E)\tag a$$
where $f_{E_v}\in\ca'$ is of the form 
$$f_{E_v}= f_0v^{-2\aa_E}+\text{ strictly higher powers of }v\tag b$$
and $f_0\in\CC-\{0\}$. (See \cite{\HEC, 19.1(e), 20.1(c), 20.7}.)

From Lemma 1.4 we see that $\l_E\i T_{w_0}$ acts on $E_v$ as $\s_E$. Using \cite{\CDGVII, 34.14(e)} with 
$c=\l_E\i T_{w_0}$ (an invertible element of $\ch_K$) we see that
$$\sum_{x\in W}\tr(T_x\s_E,E_v)\tr(\s_E\i T_{x\i},E_v)=f_{E_v}\dim(E).\tag c$$

\proclaim{Lemma 1.6} Let $E\in\Irr W$. We have $\l_E=v^{n_E}$ for some $n_E\in\ZZ$.
\endproclaim
For any $x\in W$ we have 
$$\tr(\s_Ec_x^\da,E_v)=\sum_{d\in\cd,z\in W;d=d_z}h_{x,d,z}n_d\tr(\s_Et_z,E_\spa)\in\ca'$$
since $\tr(\s_Et_z,E_\spa)\in\CC$. It follows that $\tr(\s_Eh,E_v)\in\ca'$ for any $h\in\ch$. In particular,
both $\tr(\s_ET_{w_0},E_v)$ and $\tr(T_{w_0}\i\s_E,E_v)$ belong to $\ca'$. Thus $\l_E\dim E$ and 
$\l_E\i\dim E$ belong to $\ca'$ so that $\l_E=bv^n$ for some $b\in\CC-\{0\}$ and $n\in\ZZ$. From the 
definitions we have $\l_E|_{v=1}=1$ (for $v=1$, $T_{w_0}$ becomes $w_0$) hence $b=1$. The lemma is proved.

\proclaim{Lemma 1.7} Let $E\in\Irr W$. There exists $\e_E\in\{1,-1\}$ such that for any $x\in W$ we have
$$\tr(\s_{E^\da}T_x,(E^\da)_v)=\e_E(-1)^{l(x)}\tr(\s_ET_{x\i}\i,E_v).\tag a$$
\endproclaim
Let $(E_v)^\da$ be the $\ch_K$-module with underlying vector space $E_v$ such that the action of $h\in\ch_K$ 
on $(E_v)^\da$ is the same as the action of $h^\da$ on $E_v$. From the proof in \cite{\HEC, 20.9} we see 
that there exists an isomorphism of $\ch_K$-modules $b:(E_v)^\da@>\si>>(E^\da)_v$. Let 
$\io:(E_v)^\da@>>>(E_v)^\da$ be the vector space isomorphism which  corresponds under $b$ to 
$\s_{E^\da}:(E^\da)_v@>>>(E^\da)_v$. Then we have 
$\tr(\s_{E^\da}T_x,(E^\da)_v)=\tr(\io T_x,(E_v)^\da)$.
It is enough to prove that $\io=\pm\s_E$ as a $K$-linear map of the vector space $E_v=(E_v)^\da$ into itself.
From the definition we have $\io(T_we)=T_{\s(w)}\io(e)$ for all $w\in W,e\in(E_v)^\da$. Hence 
$(-1)^{l(w)}\io(T_{w\i}\i e)=(-1)^{l(w)}T_{\s(w\i)}\i\io(e)$ for all $w\in W,e\in E_v$. It follows that 
$\io(h e)=(-1)^{l(w)}T_{\s(h)}\io(e)$ for all $h\in\ch,e\in E_v$. Hence $\io(T_we)=T_{\s(w)}\io(e)$ 
for all $w\in W,e\in E_v$. By the uniqueness in Lemma 1.4 we see that $\io=\e_E\s_E:E_v@>>>E_v$ where 
$\e_E\in K-\{0\}$. Since $\io^2=1$, $\s_E^2=1$, we see that $\e_E=\pm1$. The lemma is proved.

\proclaim{Lemma 1.8} Let $E\in\Irr W$. We have $n_E=-\aa_E+\aa_{E^\da}$. 
\endproclaim
For $x\in W$ we have (using Lemma 1.4, 1.6)
$$\tr(T_{w_0x},E_v)=\tr(T_{w_0}T_{x\i}\i,E_v)=v^{n_E}\tr(\s_ET_{x\i}\i,E_v).\tag a$$
Making a change of variable $x\m w_0x$ in 1.5(a) and using that $T_{x\i w_0}=T_{w_0\s(x)\i}$ we obtain
$$\align&f_{E_v}\dim(E)=\sum_{x\in W}\tr(T_{w_0x},E_v)\tr(T_{w_0\s(x)\i},E_v)\\&=
v^{2n_E}\sum_{x\in W}\tr(\s_ET_{x\i}\i,E_v)\tr(\s_ET_{\s(x)}\i,E_v).\endalign$$
Using now Lemma 1.7 and the equality $l(x)=l(\s(x\i))$ we obtain
$$\align&f_{E_v}\dim(E)=v^{2n_E}\sum_{x\in W}\tr(\s_{E^\da}T_x,(E^\da)_v)
\tr(\s_{E^\da}T_{\s(x\i)},(E^\da)_v)\\&=v^{2n_E}\sum_{x\in W}\tr(\s_{E^\da}T_x,(E^\da)_v)
\tr(T_{\x\i}\s_{E^\da},(E^\da)_v)\\&=v^{2n_E}f_{(E^\da)_v}\dim(E^\da).\endalign$$
(The last step uses 1.5(c) for $E^\da$ instead of $E$.) Thus we have $f_{E_v}=v^{2n_E}f_{(E^\da)_v}$.
The left hand side is as in 1.5(b) and similarly the right hand side of the form
$$f'_0v^{2n_E-2\aa_{E^\da}}+\text{strictly higher powers of }v$$
where $f_0,f'_0\in\CC-\{0\}$. It follows that $-2\aa_E=2n_E-2\aa_{E^\da}$. The lemma is proved.

\proclaim{Lemma 1.9}Let $E\in\Irr W$ and let $x\in W$. We have
$$\tr(T_x,E_v)=(-1)^{l(x)}v^{-\aa_E}\tr(t_x,E_\spa)\mod v^{-\aa_E+1}\CC[v],\tag a$$
$$\tr(\s_ET_x,E_v)=(-1)^{l(x)}v^{-\aa_E}\tr(\s_Et_x,E_\spa)\mod v^{-\aa_E+1}\CC[v].\tag b$$
\endproclaim
For a proof of (a), see \cite{\HEC, 20.6(b)}. We now give a proof of (b) along the same lines as that of (a).
There is a unique two sided cell $\boc$ such that $t_z|_{E_\spa}=0$ for $z\in W-\boc$. Let $a=\aa(z)$ for all
$z\in\boc$. By \cite{\HEC, 20.6(c)} we have $a=\aa_E$.
From the definition of $c_x$ we see that $T_x=\sum_{y\in W}f_yc_y$ where $f_x=1$ and $f_y\in v\i\ZZ[v\i]$
for $y\ne x$. Applying ${}^\da$ we obtain $(-1)^{l(x)}T_{x\i}\i=\sum_{y\in W}f_yc_y^\da$; applying $\bar{}$ 
we obtain $(-1)^{l(x)}T_x=\sum_{y\in W}\bar f_yc_y^\da$. Thus we have 
$$\align&(-1)^{l(x)}\tr(\s_ET_x,E_v)=\sum_{y\in W}\bar f_y\tr(\s_Ec_y^\da,E_v)\\&=
\sum_{y,z\in W,d\in\cd;d=d_z}\bar f_yh_{y,d,z}n_d\tr(\s_Et_z,E_\spa).\endalign$$
In the last sum we can assume that $z\in\boc$ and $d\in\boc$ so that
$h_{y,d,z}=\g_{y,d,z\i}v^{-a}\mod v^{-a+1}\ZZ[v]$. Since $\bar f_x=1$ and $\bar f_y\in v\ZZ[v]$ for all 
$y\ne x$ we see that 
$$(-1)^{l(x)}\tr(\s_ET_x,E_v)=\sum_{z\in\boc,d\in\cd\cap\boc}\g_{x,d,z\i}n_dv^{-a}\tr(\s_Et_z,E_\spa)
\mod v^{-a+1}\CC[v].$$
If $x\n\boc$ then $\g_{x,d,z\i}=0$ for all $d,z$ in the sum so that
$\tr(\s_ET_x,E_v)=0$; we have also $\tr(\s_Et_x,E_\spa)=0$ and the desired formula follows.
We now assume that $x\in\boc$. Then for $d,z$ as above we have $\g_{x,d,z\i}=0$ unless $x=z$ and $d=d_x$
in which case $\g_{x,d,z\i}n_d=1$. Thus (b) holds again. The lemma is proved.

\proclaim{Lemma 1.10} Let $E\in\Irr W$. Let $\boc$ be the unique two sided cell such that $t_z|_{E_\spa}=0$ 
for $z\in W-\boc$. Let $\boc'$ be the unique two sided cell such that $t_z|_{(E^\da)_\spa}=0$ for 
$z\in W-\boc'$. We have $\boc'=w_0\boc$.
\endproclaim
Using 1.8(a) and 1.7(a) we have
$$\tr(T_{w_0x},E_v)=v^{n_E}\tr(\s_ET_{x\i}\i,E_v)=v^{n_E}\e_E(-1)^{l(x)}\tr(\s_{E^\da}T_x,(E^\da)_v).\tag a$$
Using 1.9(a) for $E$ and 1.9(b) for $E^\da$ we obtain
$$\tr(T_{w_0x},E_v)=(-1)^{l(w_0x)}v^{-\aa_E}\tr(t_{w_0x},E_\spa)\mod v^{-\aa_E+1}\CC[v],$$
$$\tr(\s_{E^\da}T_x,(E^\da)_v)=(-1)^{l(x)}v^{-\aa_{E^\da}}\tr(\s_{E^\da}t_x,E^\da_\spa)
\mod v^{-\aa_{E^\da}+1}\CC[v].$$
Combining with (a) we obtain
$$\align&(-1)^{l(w_0x)}v^{-\aa_E}\tr(t_{w_0x},E_\spa)+\text{strictly higher powers of }v\\&=
v^{n_E}\e_Ev^{-\aa_{E^\da}}\tr(\s_{E^\da}t_x,E^\da_\spa)+\text{strictly higher powers of }v.\endalign$$
Using the equality $n_E=-\aa_E+\aa_{E^\da}$ (see 1.8) we deduce
$$(-1)^{l(w_0x)}\tr(t_{w_0x},E_\spa)=\e_E\tr(\s_{E^\da}t_x,E^\da_\spa).$$
Now we can find $x\in W$ such that $\tr(t_{w_0x},E_\spa)\ne0$ and the previous equality shows that
$t_x|_{(E^\da)_\spa}\ne0$. Moreover from the definition we have $w_0x\in\boc$ and $x\in\boc'$ so that
$w_0\boc\cap\boc'\ne\emp$. Since $w_0\boc$ is a two-sided cell (see \cite{\HEC, 11.7(d)}) it follows that 
$w_0\boc=\boc'$. The lemma is proved.

\proclaim{Lemma 1.11} Let $\boc$ be a two-sided cell of $W$. Let $\boc'$ be the two-sided cell 
$w_0\boc=\boc w_0$ (see Lemma 1.2). Let $a=\aa(x)$ for any $x\in\boc$; let $a'=\aa(x')$ for any $x'\in\boc'$.
The $K$-linear map $J^\boc_K@>>>J^\boc_K$ given by $\x\m\ph(v^{a-a'}T_{w_0})\x$ (left multiplication in
$J_K$) is obtained from a $\CC$-linear map $J^\boc_\CC@>>>J^\boc_\CC$ (with square $1$) by extension of 
scalars from $\CC$ to $K$.
\endproclaim
We can find a direct sum decomposition $J^\boc_\CC=\op_{i=1}^mE^i$ where $E^i$ are simple left ideals of 
$J_\CC$ contained in $J^\boc_\CC$. We have $J^\boc_K=\op_{i=1}^mK\ot E^i$. It is enough to show that for any
$i$, the $K$-linear map $K\ot E^i@>>>K\ot E^i$ given by the action of $\ph(v^{a-a'}T_{w_0})$ in the left 
$J_K$-module structure of $K\ot E^i$ is obtained from a $\CC$-linear map $E^i@>>>E^i$ (with square $1$) by 
extension of scalars from $\CC$ to $K$. We can find $E\in\Irr W$ such that $E^i$ is isomorphic to $E_\spa$ 
as a $J_\CC$-module. It is then enough to show that the action of $v^{a-a'}T_{w_0}$ in the left 
$\ch_K$-module structure of $E_v$ is obtained from the map $\s_E:E@>>>E$ by extension of scalars from $\CC$ 
to $K$. This follows from the equality $v^{a-a'}T_{w_0}=\s_E:E_v@>>>E_v$ (since $\s_E$ is obtained by 
extension of scalars from a $\CC$-linear map $E@>>>E$ with square $1$) provided that we show that 
$-n_E=a-a'$. Since $n_E=-\aa_E+\aa_{E^\da}$ (see Lemma 1.8) it is enough to show that $a=\aa_E$ and 
$a'=\aa_{E^\da}$. The equality $a=\aa_E$ follows from \cite{\HEC, 20.6(c)}. The equality $a'=\aa_{E^\da}$ 
also follows from \cite{\HEC, 20.6(c)} applied to $E^\da,\boc'=w_0\boc$ instead of $E,\boc$ (see Lemma 
1.10). The lemma is proved.

\proclaim{Lemma 1.12} In the setup of Lemma 1.11 we have for any $x\in\boc$:
$$\ph(v^{a-a'}T_{w_0})t_x=\sum_{x'\in\boc}m_{x',x}t_{x'}\tag a$$
$$\ph(v^{2a-2a'}T_{w_0}^2)t_x=t_x\tag b$$
where $m_{x',x}\in\ZZ$.
\endproclaim
Now (b) and the fact that (a) holds with $m_{x',x}\in\CC$ is just a restatement of Lemma 1.11. Since 
$\ph(v^{a-a'}T_{w_0})\in J_\ca$ we have also $m_{x',x}\in\ca$. We now use that $\ca\cap\CC=\ZZ$ and the 
lemma follows.

\proclaim{Lemma 1.13} In the setup of Lemma 1.11 we have for any $x\in\boc$ the following equalities in 
$\ch^\boc$:
$$v^{a-a'}T_{w_0}c_x^\da=\sum_{x'\in\boc}m_{x',x}c_{x'}^\da,\tag a$$
$$v^{2a-2a'}T_{w_0}^2c_x^\da=c_x^\da\tag b$$
where $m_{x',x}\in\ZZ$ are the same as in Lemma 1.12. Moreover, if $m_{x',x}\ne0$ then $x'\si_\cl x$.
\endproclaim
The first sentence follows from Lemma 1.12 using \cite{\HEC, 18.10(a)}.
Clearly, if $m_{x',x}\ne0$ then $x'\le_\cl x$ which together with $x'\si_{\cl\car}x$ implies $x'\si_\cl x$.

\head 2. The main results\endhead
\subhead 2.1\endsubhead
In this section we fix a two-sided cell $\boc$ of $W$; $a,a'$ are as in 1.11. We define an $\ca$-linear map 
$\th:\ch^{\le\boc}@>>>\ca$ by $\th(c_x^\da)=1$ if $x\in\cd\cap\boc$, $\th(c_x^\da)=0$ if 
$x\le_{\cl\car}x'\text{ for some }x'\in\boc$ and $x\n\cd\cap\boc$. Note that $\th$ is zero on $\ch^{<\boc}$ 
hence it can be viewed as an $\ca$-linear map $\ch^\boc@>>>\ca$.

\proclaim{Lemma 2.2} Let $x,x'\in\boc$. We have
$$\th(c_{x\i}^\da c_{x'}^\da)=n_{d_x}\d_{x,x'}v^a+\text{ strictly lower powers of }v.\tag a$$
\endproclaim
The left hand side of (a) is
$$\align&\sum_{d\in\cd\cap\boc}h_{x\i,x',d}=
\sum_{d\in\cd\cap\boc}\g_{x\i,x',d}v^a+\text{ strictly lower powers of }v\\&=
n_{d_x}\d_{x,x'}v^a+\text{ strictly lower powers of }v.\endalign$$
The lemma is proved.

We now state one of the main results of this paper. 
\proclaim{Theorem 2.3} There exists a unique permutation $u\m u^*$ of $\boc$ (with square $1$) such that for
any $u\in\boc$ we have
$$v^{a-a'}T_{w_0}c_u^\da=\e_uc_{u^*}^\da\mod\ch^{<\boc}\tag a$$ 
where $\e_u=\pm1$. For any $u\in\boc$ we have $\e_{u\i}=\e_u=\e_{\s(u)}=\e_{u^*}$ and 
$\s(u^*)=(\s(u))^*=((u\i)^*)\i$.
\endproclaim
Let $u\in\boc$. We set $Z=\th((v^{a-a'}T_{w_0}c_u^\da)^\flat v^{a-a'}T_{w_0}c_u^\da)$. We compute $Z$ in 
two ways, using Lemma 2.2 and Lemma 1.13. We have 
$$Z=\th(c_{u\i}^\da v^{2a-2a'}T_{w_0}^2c_u^\da)
=\th(c_{u\i}^\da c_u^\da)=n_{d_u}v^a+\text{ strictly lower powers of }v,$$
$$\align&Z=\th((\sum_{y\in\boc}m_{y,u}c_y^\da)^\flat(\sum_{y'\in\boc}m_{y',u}c_{y'}^\da))=
\sum_{y,y'\in\boc}m_{y,u}m_{y',u}\th(c_{y\i}^\da c_{y'}^\da)\\&=
\sum_{y,y'\in\boc}m_{y,u}m_{y',u}n_{d_y}\d_{y,y'}v^a+\text{ strictly lower powers of }v\\&
=\sum_{y\in\boc}n_{d_y}m_{y,u}^2v^a+\text{ strictly lower powers of }v\\&
=\sum_{y\in\boc}n_{d_u}m_{y,u}^2v^a+\text{ strictly lower powers of }v\endalign$$
where $m_{y,u}\in\ZZ$ is zero unless $y\si_\cl u$ (see 1.13), in which case we have $d_y=d_u$. We deduce 
that $\sum_{y\in\boc}m_{y,u}^2=1$, so that we have $m_{y,u}=\pm1$ for a unique $y\in\boc$ (denoted by $u^*$)
and $m_{y,u}=0$ for all $y\in\boc-\{u^*\}$.
Then (a) holds. Using (a) and Lemma 1.13(b) we see that $u\m u^*$ has square $1$ and that $\e_u\e_{u^*}=1$.

The automorphism $\s:\ch@>>>\ch$ (see 1.1) satisfies the equality $\s(c_u^\da)=c_{\s(u)}^\da$ for any 
$u\in W$; note also that $w\in\boc\lra\s(w)\in\boc$ (see Lemma 1.2). Applying $\s$ to (a) we obtain
$$v^{a-a'}T_{w_0}c_{\s(u)}^\da=\e_uc_{\s(u^*)}^\da$$ 
in $\ch^\boc$.
By (a) we have also $v^{a-a'}T_{w_0}c_{\s(u)}^\da=\e_{\s(u)}c_{(\s(u))^*}^\da$ in $\ch^\boc$.
It follows that $\e_uc_{\s(u^*)}^\da=\e_{\s(u)}c_{(\s(u))^*}^\da$ hence 
$\e_u=\e_{\s(u)}$ and $\s(u^*)=(\s(u))^*$. 

Applying $h\m h^\flat$ to (a) we obtain
$$v^{a-a'}c_{u\i}^\da T_{w_0}=\e_uc_{(u^*)\i}^\da$$ 
in $\ch^\boc$. By (a) we have also
$$v^{a-a'}c_{u\i}^\da T_{w_0}=v^{a-a'}T_{w_0}c_{\s(u\i)}^\da=\e_{\s(u\i)}c_{(\s(u\i))^*}^\da$$
in $\ch^\boc$. It follows that $\e_uc_{(u^*)\i}^\da=\e_{\s(u\i)}c_{(\s(u\i))^*}^\da$ hence 
$\e_u=\e_{\s(u\i)}$ and $(u^*)\i=(\s(u\i))^*$. Since $\e_{\s(u\i)}=\e_{u\i}$, we see that $\e_u=\e_{u\i}$.
Replacing $u$ by $u\i$ in $(u^*)\i=(\s(u\i))^*$ we obtain $((u\i)^*)\i=(\s(u))^*$ as required.
The theorem is proved. 

\subhead 2.4\endsubhead
For $u\in\boc$ we have
$$u\si_\cl u^*,\tag a$$
$$\s(u)\si_\car u^*.\tag b$$
Indeed, (a) follows from 1.13. To prove (b) it is enough to show that \lb
$\s(u)\i\si_\cl(u^*)\i$. Using (a) for $\s(u)\i$ instead of $u$ we see that it is enough to show that 
$(\s(u\i))^*=(u^*)\i$; this follows from 2.3.

If we assume that 

(c) {\it any left cell in $\boc$ intersects any right cell in $\boc$ in exactly one element}
\nl
then by (a),(b), for any $u\in\boc$, 

(d) {\it $u^*$ is the unique element of $\boc$ in the intersection of the left cell of $u$ with right cell
of $\s(u)$.}
\nl
Note that condition (c) is satisfied for any $\boc$ if $W$ is of type $A_n$ or if $W$ is of type $B_n$
($n\ge2$) with $L(s)=2$ for all but one $s\in S$ and $L(s)=1$ or $3$ for the remaining $s\in S$. (In this last case we are in the quasisplit case and
we have $\s=1$ hence $u^*=u$ for all $u$.)

\proclaim{Theorem 2.5} For any $x\in W$ we set $\vt(x)=\g_{w_0d_{w_0x\i},x,(x^*)\i}$.

(a) If $d\in\cd$ and $x,y\in\boc$ satisfy $\g_{w_0d,x,y}\ne0$ then $y=(x^*)\i$.

(b) If $x\in\boc$ then there is a unique $d\in\cd\cap w_0\boc$ such that $\g_{w_0d,x,(x^*)\i}\ne0$, namely
$d=d_{w_0x\i}$. Moreover we have $\vt(x)=\pm1$.

(c) For $u\in\boc$ we have $\e_u=(-1)^{l(w_0d)}n_d\vt(u)$ where $d=d_{w_0u\i}$. 
\endproclaim
Appplying $h\m h^\da$ to 2.3(a) we obtain for any $u\in\boc$:
$$v^{a-a'}(-1)^{l(w_0)}\ov{T_{w_0}}c_u=\sum_{z\in\boc}\d_{z,u^*}\e_u c_z\mod\sum_{z'\in W-\boc}\ca c_{z'}.
\tag d$$
We have $T_{w_0}=\sum_{y\in W}(-1)^{l(w_0y)} p_{1,w_0y}c_y$ hence
$\ov{T_{w_0}}=\sum_{y\in W}(-1)^{l(w_0y)}\ov{p_{1,w_0y}}c_y$. Introducing this in (d) we obtain
$$v^{a-a'}\sum_{y\in W}(-1)^{l(y)}\ov{p_{1,w_0y}}c_yc_u
=\sum_{z\in\boc}\d_{z,u^*}\e_u c_z\mod\sum_{z'\in W-\boc}\ca c_{z'}$$
that is,
$$v^{a-a'}\sum_{y,z\in W}(-1)^{l(y)}\ov{p_{1,w_0y}}h_{y,u,z}c_z
=\sum_{z\in\boc}\d_{z,u^*}\e_u c_z\mod\sum_{z'\in W-\boc}\ca c_{z'}.$$
Thus, for $z\in\boc$ we have
$$v^{a-a'}\sum_{y\in W}(-1)^{l(y)}\ov{p_{1,w_0y}}h_{y,u,z}=\d_{z,u^*}\e_u.\tag e$$
Here we have $h_{y,u,z}=\g_{y,u,z\i}v^{-a}\mod v^{-a+1}\ZZ[v]$ and we can assume than $z\le_{\car}y$
so that $w_0y\le_{\car}w_0z$ and $\aa(w_0y)\ge\aa(w_0z)=a'$.

For $w\in W$ we set $s_w=n_w$ if $w\in\cd$ and $s_w=0$ if $w\n\cd$. By \cite{\HEC, 14.1} we have
$p_{1,w}=s_wv^{-\aa(w)}\mod v^{-\aa(w)-1}\ZZ[v\i]$ hence $\ov{p_{1,w}}=s_wv^{\aa(w)}\mod v^{\aa(w)+1}\ZZ[v]$.
Hence for $y$ in the sum above we have $\ov{p_{1,w_0y}}=s_{w_0y}v^{\aa(w_0y)}\mod v^{\aa(w_0y)+1}\ZZ[v]$.
Thus (e) gives
$$v^{a-a'}\sum_{y\in\boc}(-1)^{l(y)}s_{w_0y}\g_{y,u,z\i}v^{\aa(w_0y)-a}-\d_{z,u^*}\e_u\in v\ZZ[v]$$
and using $\aa(w_0y)=a'$ for $y\in\boc$ we obtain
$$\sum_{y\in\boc}(-1)^{l(y)}s_{w_0y}\g_{y,u,z\i}=\d_{z,u^*}\e_u.$$
Using the definition of $s_{w_0y}$ we obtain
$$\sum_{d\in\cd\cap w_0\boc}(-1)^{l(w_0d)}n_d\g_{w_0d,u,z\i}=\d_{z,u^*}\e_u.\tag f$$
Next we note that 

(g) {\it if $d\in\cd$ and $x,y\in\boc$ satisfy $\g_{w_0d,x,y}\ne0$ then $d=d_{w_0 x\i}$.}
\nl
Indeed from \cite{\HEC,\S14, P8} we deduce $w_0d\si_\cl x\i$. Using \cite{\HEC, 11.7} we deduce
$d\si_\cl w_0 x\i$ so that $d=d_{w_0\i x\i}$. This proves (g).

Using (g) we can rewrite (f) as follows.
$$(-1)^{l(w_0)}(-1)^{l(d)}n_d\g_{w_0d,u,z\i}=\d_{z,u^*}\e_u\tag h$$
where $d=d_{w_0u\i}$.

We prove (a). Assume that $d\in\cd$ and $x,y\in\boc$ satisfy $\g_{w_0d,x,y}\ne0$, $y\ne(x^*)\i$.
Using (g) we have $d=d_{w_0 x\i}$. Using (h) with $u=x,z=y\i$ we see that $\g_{w_0d,x,y}=0$, a contradiction.
This proves (a).

We prove (b). Using (h) with $u=x,z=x^*$ we see that
$$(-1)^{l(w_0d)}n_d\g_{w_0d,x,(x^*)\i}=\e_u\tag i$$
where $d=d_{w_0x\i}$.
Hence the existence of $d$ in (b) and the
equality $\vt(x)=\pm1$ follow; the uniqueness of $d$ follows from (g).

Now (c) follows from (i). This completes the proof of the theorem.

\subhead 2.6\endsubhead
In the case where $L=l$, $\vt(u)$ (in 2.5(c)) is $\ge0$ and $\pm1$ hence $1$; moreover, 
$n_d=1$, $(-1)^{l(d)}=(-1)^{a'}$ for any $d\in\cd\cap w_0\boc$ (by the definition of $\cd$).
Hence we have $\e_u=(-1)^{l(w_0)+a'}$ for any $u\in\boc$, a result of \cite{\MAT}.

Now Theorem 2.5 also gives a characterization of $u^*$ for $u\in\boc$; it is the unique element $u'\in\boc$ 
such that $\g_{w_0d,u,u'{}\i}\ne0$ for some $d\in\cd\cap w_0\boc$.

We will show:

(a) {\it The subsets $X=\{d^*;d\in\cd\cap\boc\}$ and $X'=\{w_0d';d'\in\cd\cap w_0\boc\}$ of $\boc$ coincide.}
\nl
Let $d\in\cd\cap\boc$. By 2.5(b) we have $\g_{w_0d',d,(d^*)\i}=\pm1$ for some $d'\in\cd\cap w_0\boc$. Hence 
$\g_{(d^*)\i,w_0d',d}=\pm1$. Using \cite{\HEC, 14.2, P2} we deduce $d^*=w_0d'$. Thus $X\sub X'$. Let $Y$ 
(resp. $Y'$) be the set of left cells contained in $\boc$ (resp. $w_0\boc$). We have $\sha(X)=\sha(Y)$ and 
$\sha(X')=\sha(Y')$. By \cite{\HEC, 11.7(c)} we have $\sha(Y)=\sha(Y')$. It follows that $\sha(X)=\sha(X')$.
Since $X\sub X'$, we must have $X=X'$. This proves (a).

\proclaim{Theorem 2.7} We have
$$\ph(v^{a-a'}T_{w_0})=\sum_{d\in\cd\cap\boc}\vt(d)\e_dt_{d^*}\mod\sum_{u\in W-\boc}\ca t_u.$$
\endproclaim
We set $\ph(v^{a-a'}T_{w_0})=\sum_{u\in W}p_ut_u$ where $p_u\in\ca$.
Combining 1.12(a), 1.13(a), 2.3(a) we see that for any $x\in\boc$ we have 
$$\ph(v^{a-a'}T_{w_0})t_x=\e_xt_{x^*},$$
hence
$$\e_xt_{x^*}=\sum_{u\in\boc}p_ut_ut_x=\sum_{u,y\in\boc}p_u\g_{u,x,y\i}t_y.$$
It follows that for any $x,y\in\boc$ we have
$$\sum_{u\in\boc}p_u\g_{u,x,y\i}=\d_{y,x^*}\e_x.$$
Taking $x=w_0d$ where $d=d_{w_0y}\in\cd\cap w_0\boc$ we obtain 
$$\sum_{u\in\boc}p_u\g_{w_0d_{w_0y},y\i,u}=\d_{y,(w_0d_{w_0y})^*}\e_{w_0d_{w_0y}}$$
which, by 2.5, can be rewritten as 
$$p_{((y\i)^*)\i}\vt(y\i)=\d_{y,(w_0d_{w_0y})^*}\e_{w_0d_{w_0y}}.$$
We see that for any $y\in\boc$ we have
$$p_{\s(y^*)}=\d_{y,(w_0d_{w_0y})^*}\vt(y\i)\e_{w_0d_{w_0y}}.$$
In particular we have $p_{\s(y^*)}=0$ unless $y=(w_0d_{w_0y})^*$ in which case
$$p_{\s(y^*)}=p_{(\s(y))^*)}=\vt(y\i)\e_y.$$
(We use that $\e_{y^*}=\e_y$.) If $y=(w_0d_{w_0y})^*$ then
$y^*\in X'$ hence by 2.6(a), $y^*=d^*$ that is $y=d$ for some $d\in\cd$. Conversely, if $y\in\cd$ then
$w_0y^*\in\cd$ (by 2.6(a)) and $w_0y^*\si_\cl w_0y$ (since $y^*\si_\cl y$) hence $d_{w_0y}=w_0y^*$. We see
that $y=(w_0d_{w_0y})^*$ if and only if $y\in\cd$. We see that
$$\ph(v^{a-a'}T_{w_0})=\sum_{d\in\cd\cap\boc}\vt(d\i)\e_dt_{(\s(d))^*}+\sum_{u\in W-\boc}p_ut_u.$$
Now $d\m\s(d)$ is a permutation of $\cd\cap\boc$ and $\vt(d\i)=\vt(d)=\vt(\s(d))$, $\e_{\s(d)}=\e_d$. 
The theorem follows.

\proclaim{Corollary 2.8} We have
$$\ph(T_{w_0})=\sum_{d\in\cd}\vt(d)\e_dv^{-\aa(d)+\aa(w_0d)}t_{d^*}\in J_\ca.$$
\endproclaim

\subhead 2.9\endsubhead
We set $\fT_\boc=\sum_{d\in\cd\cap\boc}\vt(d)\e_dt_{d^*}\in J^\boc$. We show:

(a) $\fT_\boc^2=\sum_{d\in\cd\cap\boc}n_dt_d$;

(b) $t_x\fT_\boc=\fT_\boc t_{\s(x)}$ for any $x\in W$.
\nl
By 2.7 we have $\ph(v^{a-a'}T_{w_0})=\fT_\boc+\x$ where $\x\in J^{W-\boc}_K:=\sum_{u\in W-\boc}Kt_u$. Since 
$J^\boc_K, J^{W-\boc}_K$ are two-sided ideals of $J_K$ with intersection zero
and $\ph_K:\ch_K@>>>J_K$ is an algebra homomorphism, it follows that
$$\ph(v^{2a-2a'}T_{w_0}^2)=(\ph(v^{a-a'}T_{w_0}))^2=(\fT_\boc+\x)^2=\fT^2_\boc+\x'$$
where $\x'\in J^{W-\boc}_K$.  Hence, for any $x\in\boc$  we have 
$\ph(v^{2a-2a'}T_{w_0}^2)t_x=\fT^2_\boc t_x$ so that (using 1.12(b)):
$t_x=\fT^2_\boc t_x$. We see that $\fT^2_\boc$ is the unit element of the ring $J^\boc_K$. Thus (a) holds.

We prove (b). For any $y\in W$ we have $T_yT_{w_0}=T_{w_0}T_{\s(y)}$ hence, applying $\ph_K$,
$$\ph(T_y)\ph(v^{a-a'}T_{w_0})=\ph(v^{a-a'}T_{w_0})\ph(T_{\s(y)})$$
that is, $\ph(T_y)(\fT_\boc+\x)=(\fT_\boc+\x)\ph(T_{\s(y)})$. Thus,
$\ph(T_y)\fT_\boc=\fT_\boc\ph(T_{\s(y)})+\x_1$ where $\x_1\in J^{W-\boc}_K$.
Since $\ph_K$ is an isomorphism, it follows that for any $x\in W$ we have 
$t_x\fT_\boc=\fT_\boc t_{\s(x)}\mod J^{W-\boc}_K$. Thus (b) holds.

\subhead 2.10\endsubhead
In this subsection we assume that $L=l$. In this case 2.8 becomes
$$\ph(T_{w_0})=\sum_{d\in\cd}(-1)^{l(w_0)+\aa(w_0d)}v^{-\aa(d)+\aa(w_0d)}t_{d^*}\in J_\ca.$$
(We use that $\vt(d)=1$.)

\mpb

For any left cell $\G$ contained in $\boc$ let $n_\G$ be the number of fixed 
points of the permutation $u\m u^*$ of $\G$. Now $\G$ carries a representation $[\G]$ of $W$ and from 2.3 we
see that $\tr(w_0,[\G])=\pm n_\G$. Thus $n_\G$ is the absolute value of the integer $\tr(w_0,[\G])$. From 
this the number $n_\G$ can be computed for any $\G$. In this way we see for example that if $W$ is of type 
$E_7$ or $E_8$ and $\boc$ is not an exceptional two-sided cell, then $n_\G>0$.

\enddocument